\pdfoutput=1
\documentclass{article}
\usepackage{graphicx} 

\usepackage{graphicx} 
\usepackage{amssymb}  
\usepackage{amsmath}
\newcommand{\ol}{\overline}

\newcommand{\Pow}{{\cal P}}
\newcommand{\s}{\subseteq}
\newcommand{\es}{\emptyset}
\newcommand{\Ra}{\Rightarrow}
\newcommand{\LRa}{\Leftrightarrow}

\title{Compression with wildcards: All metric induced subgraphs}

\author{Marcel Wild}

\begin{document}

\maketitle

\begin{quote}
A{\scriptsize BSTRACT}: {\footnotesize Driven by applications in the natural, social and computer sciences several algorithms have been proposed to enumerate all sets $X\s V$ of vertices of a graph $G=(V,E)$ that induce a {\it connected} subgraph. We offer two algorithms for  
enumerating all $X$'s that induce (more exquisite) {\it metric} subgraphs.  Specifically, the first algorithm, 
called {\tt AllMetricSets}, generates these $X$'s in a compressed format.
The second algorithm generates all (accessible) metric sets one-by-one but is provably output-polynomial.  Mutatis mutandis the same holds for the geodesically convex
sets $X\s V$; this being  a natural strengthening of "metric". The Mathematica command {\tt BooleanConvert} features prominently.
  }
\end{quote}

\section{Introduction}
All graphs $G=(V,E)$ will be undirected and simple. Here $V$ is the vertex set and $E$ the edge set.
Consider the graph $G$ in Figure 1 below. By abuse of language the subset of vertices 
$X:=\{1,2,4,5,6,7,8\}\s V:=\{1,2,...,12\}$ is {\it connected} in the sense that the subgraph induced\footnote{Formally the subgraph  {\it induced} by a set $X\s V$ is the graph $G[X]:=(X,E')$ with edge set  $E':=\{\{i,j\}\in E:\ i,j\in X\}$. For our $X$  the induced edges are rendered boldface in Figure 1. }  by $X$ is connected in the usual sense. 
However, $X$ in Fig.1 is {\it not metric}\footnote{Formal definitions of "metric" and "distance" will be given in Subsection 2.1.} in the sense that the distance between $1$ and $8$  within $X$ is larger than within the whole graph; indeed the path (1,5,6,7,8) (also (1,2,6,7,8)) is longer than (1,2,3,8).

The author's previous research on Boolean functions  facilitated the  design of the algorithm 
{\tt AllMetricSets} which enumerates all metric subsets of any graph. The initial motivation  was  the previous enumeration of (ordinarily) connected subsets by various authors. One perk of {\tt AllMetricSets} is the compressed enumeration, i.e. possibly millions of metric subsets $X$ are "bundled", i.e. not output one-by-one. Nonetheless each individual $X$ is readily available, if desired. 

\vspace{2mm}
\includegraphics[scale=0.85]{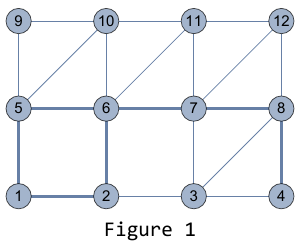}
\vspace{2mm}

 \vspace{2mm}
 Here comes the Section break-up.
In Section 2 some literature concerning the relevance of connected sets (of vertices) is surveyed. In Section 3 we sketch why the family $GeConv(G)$ of all geodesically-convex sets can be enumerated in output-polynomial time.
The familily $Met(G)$ of all metric sets of $G$ is more stubborn in this regard. However the subfamily $AccMet(G)\s Met(G)$ of all "accessible" metric sets can be enumerated in output-polynomial time as well.
Furthermore we pinpoint scenarios where $AccMet(G)= Met(G)$.

The remaining Sections are dedicated to the {\it compressed}, i.e. chunk-wise, enumeration of $Met(G)$ and $GeConv(G)$. It is natural to start (Section 4) with trees $G=T$, because in this scenario all three notions (connected, metric, geodesically-convex) mean the same.
In Section 5 we introduce our algorithm {\tt AllMetricSets}.
While it lacks output-polynomiality, it excels with  compressed enumeration.
 Section 6 numerically evaluates high-level Mathematica implementations of {\tt AllMetricSets} and {\tt AllGeConvexSets}, and compares them with a straightforward application of the hardwired Mathematica command {\tt BooleanConvert}.

\section{Community detection}

The issue of "community detection" in networks (=graphs) has become a prominent research theme in the new millenium, as larger and larger datasets have actually been implemented as graphs. We recommend [FH] as a survey of the many ways (similarity measures, probability measures, clique-related concepts, etc) that have been tried to capture the concept of "community".

 Each clique $X\s V$ is of course the most extreme (and hence unlikely) kind of connected\footnote{Note that the concept of a "connected set of vertices" does not yet occur in [FH].} set.  Let us compare the two a bit more. In contrast to connected sets, each subset of a clique is a clique, and so one may be content to enumerate only the {\it maximal cliques}. This issue is well understood and has a long history. 
 
 Concerning the connected sets of $G$, the state of the art [U]  is that they can be enumerated 
 with constant delay $O(1)$ (between any two output sets). While this is better than the linear delay algorithm of [ASA], the latter reference refines "connected" to "connected {\it and} cohesive", and also provides some bio-science background.

\vspace{2mm}
It becomes more complicated when only the sets of fixed cardinality $k$ are sought; we recommend  [WX] for several reasons: the brief summaries of previous articles on that matter; discussing published applications in bioinformatics and elsewhere; the use of toy examples; the comprehensive computer experiments.

 \vspace{2mm}
 {\bf 2.1} There are applications where  $X\s V$ being connected is not enough, i.e. it can be important to get from any $s\in X$ to any $t\in X$ fast (meaning there is no benefit in momentarily leaving $X$).

Formally one defines $dist_X(s,t)$ as the length $k$ of a shortest $s-t$ path $(v_0,v_1,...,v_k)$ (with $s=v_0,\ t=v_k$) that only uses vertices  $v_i\in X$. Evidently $dist_X(s,t)\ge dist_V(s,t)$ for all $s,t\in X$. Following [H,p.324] call
 $Y\s V$  {\it metric} if $dist_Y=dist_V$. Put another way, for each $s\neq t$ in $Y$ there must be at least one (globally) shortest $s-t$ path that lies\footnote{By saying that some $s-t$ path "lies in" $Y$ we henceforth mean that all its vertices $s,v_1,..,t$ lie in $Y$.} in $Y$.
  Following [H,p.324] we henceforth call a shortest $s-t$ path a $s-t$ {\it geodesic}. Even better than metric sets are {\it geodesically\footnote{One cannot simply say "convex" since there are monophonically-convex sets as well [W3].}-convex} sets $Y\s V$, i.e. by definition $Y$ contains, for all distinct $s,t\in Y$, {\it all} $s-t$ geodesics.

\section{Topological properties of $Met(G)$ and $Conv(G)$}

Let $Met(G)$ and $GeConv(G)$ be the families of all sets $X\s V$ that induce a metric, respectively geodesically-convex, subgraph $G[X]$. In 3.2 (having had a closer look at geodesics in 3.1) we are going to show that $GeConv(G)$ can be enumerated in output-polynomial time. The matter for $Met(G)$ is more cumbersome and will be tackled in 3.3 and 3.5. Subsection 3.4 relates to all of that and in its core is about distance-hereditary and Ptolemaic graphs.

\vspace{2mm}
{\bf 3.1} We defined metric vertex subsets of $V$  in terms of geodesics. Conversely the underlying vertex set of any geodesic is itself\footnote{To spell this out, consider any geodesic $P:=(s=v_0,v_1,\cdots,v_k=t)$. Let $i<j$. Assuming there was a $v_i-v_j$ path shorter than the subpath $(v_i,\ldots,v_j)$ of $P$, leads to the contradiction that there is some $s-t$ path shorter than our $s-t$ geodesic $P$.} a metric set. Fortunately enumerating {\it these}  metric sets has a long history, see [S]. A short survey that focuses on exploiting the distance matrix $D(G)$ (and less on how to get it) can be found in [W2]. Article [W2] also features a novel algorithm for calculating {\it all} $N'$ geodesics of $G=(V,E)$. Like the old method, it runs in output-polynomial time $O(N'|V|^3)$, but is more elegant, and  the more competitive the sparser $G$. 

Whatever method generates all metric sets of $G$, generates "along the way" all geodesics of $G$. Our particular method {\tt AllMetricSets} demands that all geodesics be generated first, because they are essential to find the remainder of $Met(G)$.

\vspace{3mm}
{\bf 3.2}   It is easy to see that  $GeConv(G)$ is a {\it closure system}, i.e. from $X,Y\in GeConv(G)$ follows $X\cap Y\in GeConv(G)$. Better still, it is easy to pinpoint an "implication-base" for $GeConv(G)$. Whenever  an explicite implication-base for a closure system $\cal C$ is known, one can enumerate $\cal C$ in ouput-polynomial time (and additionally in a compressed format). In our present scenario this yields

\vspace{2mm}
{\bf Theorem 1: }{\it If $N:=|GeConv(G)|$, then $GeConv(G)$ can be enumerated in output-polynomial time $O(Nn^6)$.}

\vspace{2mm} For a proof of Theorem 1 and background on implication-bases see [W3]. The underlying algorithm {\tt AllGeConvSets} will undergo numerical experiments in Section 6.

\vspace{3mm}
{\bf 3.3} The "topology" of $Met(G)$ is more intricate than the topology of\\ $GeConv(G)\s Met(G)$ because $Met(G)$ is no closure system, and neither does it follow from ($X\in Met(G)$ and $Y\s X$)  that $Y\in Met(G)$. 

In order to get leverage for proving an enumeration result, we call a metric set $X\s V$ {\it accessible} if there is a tower of metric sets $\es\subset Y_1\subset Y_2\subset\cdots\subset Y_k=X$ such that $|Y_i|=i$ for all $1\le i\le k$. We let $AccMet(G)\s Met(G)$ be the family of all acc-metric (= accessible metric) sets of $G$. Examples of graphs $G$ where $"\s"$ becomes $"="$ will be provided in 3.5.
 
\vspace{2mm}
{\bf Theorem 2: }{\it If $G=(V,E)$ and $N:=|AccMet(G)|$, then $AccMet(G)$ can be enumerated in output-polynomial time $O(N^2n^3)$.}

\vspace{2mm}
{\it Proof.}  Each singleton $\{x\}\ (x\in V)$ is an acc-metric set. By induction assume that all $k$-element acc-metric sets $X_1,...,X_{p(k)}$ are known. By definition of accessibility each $(k+1)$-element acc-metric set is an extension of some $X_i$. Therefore
for each fixed $X_i$ and all $t\in V\setminus X_i$ we initialize a "book-keeping" set  as $B_t:=\es$ and inflate it as follows. For each of the $N'$ geodesics $P$ (i.e. their underlying vertex sets) check whether $P=(P\cap X_i)\cup\{t\}$. If yes, replace\footnote{Hence for each $s\in B_t\s X_i$ there is a $s-t$ geodesic that lies in $X_i\cup\{t\}$.} $B_t$ by $B_t\cup (P\cap X_i)$.  Having scanned all $P$'s, check whether $B_t=X_i$. If yes then $X_i\cup\{t\}$ is acc-metric, otherwise it isn't. In the first case add $X_i\cup\{t\}$ to the list of acc-metric $(k+1)$-element sets found so far.

As to the cost analysis, first note that an acc-metric $(k+1)$-element set $X_i\cup\{y\}$ may have occured beforehand in the guise of $X_j\cup\{z\}$ for some $j<i$. Yet we may not list it twice (if only for the sake of the present proof). Fortunately double listing can be avoided as follows. Putting

$$Avoid(i):=\Bigl\{y\in V:\ (\exists j<i)\ (X_j\setminus X_i)=\{y\}\Bigr\}$$

\noindent
one checks that for each fixed $y\in V\setminus X_i$ the following holds:\\ 
($X_i\cup\{y\}=X_j\cup\{z\}$ for some $j<i$) iff $y\in Avoid(i)$. Therefore to avoid double listing it is necessary to pick $y$ in
$V\setminus (X_i\cup Avoid(i))$ Of course only the sets $X_i\cup\{y\}$ that happen to be metric are retained. 

In view of $i\le p(k)\le N$ calculating a fixed set $Avoid(i)$ costs $O(Nn)$. And deciding (by inflating $B_y$ above) whether $X_i\cup\{y\}$ is acc-metric costs $O(N'n)$.

For each $X_i\in\{X_1,X_2,...,X_{p(k)}\}$ it must be decided  whether or not $X_i$ yields at least one new acc-metric set $X_i\cup\{y\}$. Whatever the outcome, one such decision costs $O(Nn)+O(N'n)=O(Nn)$, in view of $N'\le N$.
Hence handling the whole $k$-level $\{X_1,X_2,...,X_{p(k)}\}$ costs $O(p(k)Nn)$. Since $p(1)+p(2)+\cdots+p(k_{max})=N$, the whole procedure costs is $O(N^2n)$.  Taking into account the cost $O(N'n^3)$ for generating the $N'$ geodesics, the algorithm's total cost amounts to $O(N'n^3)+O(N^2n)=O(N^2n^3)$.       $\square$

\vspace{2mm}
From the proof it is clear that outputting all $N_0$ many acc-metric sets of cardinality $\le k$ costs $O(N_0^2n^3)$. See [W3] for more scenarios of this type.

\vspace{3mm}
{\bf 3.4 }
For each graph $G=(V,E)$ and all $X\s V$ obviously these implications of properties take place:

$$(2)\quad geodesically-convex\ \Ra\ metric\ \Ra\ connected$$

None of the implications in (2) can be reversed. However, that changes when we restrict the type of graph considered. Most obviously, for trees or complete graphs all three properties in (2) are equivalent. Henceforth $Conn(G)$ denotes the family of all connected sets $X\s V$.

\vspace{3mm}
{\bf 3.4.1 }
A graph $G$ is {\it distance-hereditary} [H,p.324] if each connected set $X\s V$ is metric. In other words,  $Met(G)\s Conn(G)$ becomes
$Met(G)= Conn(G)$. One can show [H, p.325] that 

\begin{itemize}
    \item[(3)] 
{\it A graph is distance-hereditary iff each circuit of length $\ge 5$ has at least two  chords  that cross each other.}
\end{itemize}

This is illustrated by the graphs $G_2,\ G_3$ in Figures 2 and 3 respectively. The length 5 circuit of $G_2$ has  two chords, but they do not cross each other. Taking $X=\{a,b,c,d\}$ it indeed holds that $dist_X(a,d)=3>2=dist_V(a,d)$, and so $X$ is not metric. In contrast one checks that $G_3$ with its two intersecting chords is distance-hereditary.

\begin{figure}[htbp]
    \centering
    \includegraphics[width=1\linewidth]{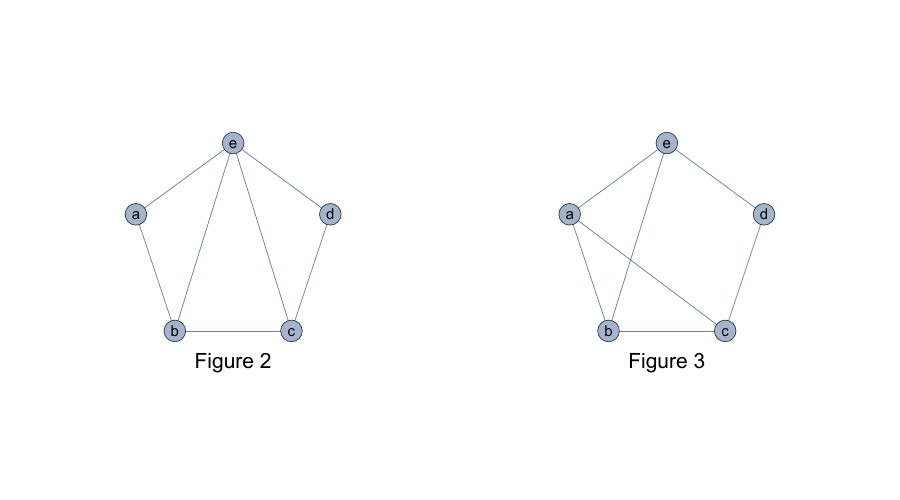}
\end{figure}

It follows that, restricted\footnote{Actually, as seen in 5.5, this restriction can easily be avoided.} to distance-hereditary graphs,  {\tt AllMetricSets} (discussed fully in Section 5) enumerates not just some "exquisite" connected sets but {\it all} connected sets; and in compressed form. Recall that all precursors, even the champion [U], enumerate the connected sets one-by-one.

\vspace{2mm}
{\bf 3.4.2} The word "convex" can mean many things in mathematics, for instance when speaking of "convex programming". Yet another one  is the concept [EJ] of a {\it convex geometry} $\cal C$. Roughly speaking, this is a  frequently occuring type of closure system with nice properties; among other things, $\cal C$  is accessible in the same sense as  $AccMet(G)$. 

It turns out that the closure system  $GeConv(G)$  is a convex geometry in the [EJ] sense iff the underlying graph is a so-called {\it Ptolemaic} graph.
The latter are defined [EJ,p.256] by an equality modelled after the  Theorem of Ptolemy (which is beautiful and e.g. subsumes Pythagoras' Theorem as an immediate  special case). Each Ptolemaic graph turns out to be distance-hereditary. In fact

$$Ptolemaic\ \LRa\ (distance\!-\!hereditary\ and\ chordal).$$

In view of (3) we can put it this way. Distance-hereditary graphs are shy of being chordal because of their misbehaving 4-circuits, whereas Ptolemaic graphs are chordal in a lush way. See  [BM],[H] for further results. 

\vspace{2mm}
{\bf 3.5} Here we further discuss acc-metric sets. For starters, not all metric sets are acc-metric; if $G=(V,E)$ is a chordless circuit with at least  5 elements then $X:=V$ is metric but not acc-metric since no set $X\setminus\{x\}$ is metric. 

There are however several sufficient\footnote{We leave it to the reader to ponder applications where $AccMet(G)$ is more desirable than $Met(G)$ (independent of whether the two coincide).} conditions for $AccMet(G)=Met(G)$. Most strikingly (but how frequently?), "=" is implied if the only metric sets are the geodesics! Likewise, if e.g. each metric set is  a tree (respectively: is complete). 
Slightly more subtle, if $Met(G)=Conn(G)$ (i.e., 3.4.1, $G$ is distance-hereditary), then $AccMet(G)=Met(G)$ takes place as well. That's because\footnote{That is best argued by considering leaves of spanning trees of connected induced subgraphs, as spelled out in [AF, Lemma 3.7].} $Conn(G)$ is accessible for whatever graph $G$. Put another way, we can say that each distance-hereditary graph enjoys this property:

\begin{itemize}
    \item[(4)]{\it For each metric $X\s V$ of the graph $G=(V,E)$ there is $x\in X$ such that $X\setminus\{x\}$ is metric.}
\end{itemize}

\noindent
Evidently {\it each} graph $G$ (distance-hereditary or not) that satisfies (4) also satisfies $AccMet(G)=Met(G)$. For instance $G_2$ in Figure 2, which is not distance-hereditary, is easily seen to satisfy
(4).

\section{Compressing $Conn(T)$ for  trees $T$}

 Apart from visual appeal\footnote{In particular that serves us well to illustrate, in Fig.4, how compression works.}, the  reason for paying special attention to trees $G=T$ is that for any two vertices $s\neq t$ in a tree there is exactly one path connecting them, which hence is the unique $s-t$ geodesic. Thus for all $Y\s V$ it holds that  

$$Y\ is\ geodesically-convex\ \LRa\ Y\ is\ connected\ \LRa\ Y\ is\footnote{ More precisely: Y induces a subtree (as opposed to just a subforest).}\ a\ subtree$$

We can hence use {\tt AllGeConvSets} to enumerate $Conn(T)=GeConv(T)$. 
\vspace{2mm}
{\bf 4.1} But back to trees.
 Consider  $T$ in Figure 4A. Lots of subtrees $Y\s V=\{1,...,10\}$ can be obtained as follows. Take the blue vertices 1,2,5,10 and add any number of white vertices. Then the resulting sets $Y$ of vertices are  subtrees. The family of these $2^6$ subtrees can be written, in obvious ways\footnote{The 2's can freely assume the values 0 or 1. Hence e.g.
 $(1,1,{\bf 0,0},1,{\bf 0,1,0,1},1)$, which matches the subtree $\{1,2,5,7,9,10\}$ is a member of $r_1$.}, as the   {\it 012-row} $r_1:=(1,1,2,2,1,2,2,2,2,1)$.

\begin{figure}[htbp]
    \centering
    \includegraphics[width=\linewidth]{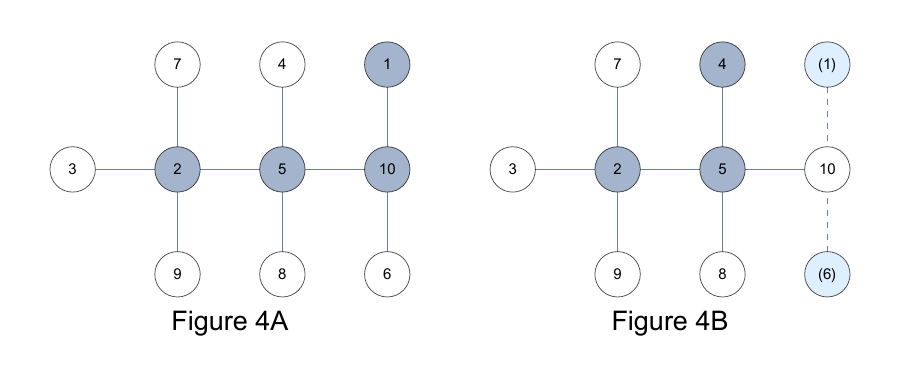}
\end{figure}

Similarly each set $Y$ in $r_2:=(0,1,2,1,1,0,2,2,2,2)$  is a subtree (see Fig.4B). Notice that $1\in Y$ for all $Y\in r_1$, while $1\not \in Y$ for all $Y\in r_2$.
Hence $r_1\cap r_2=\es$, and so we found already $64+32=96$ subtrees of $T$. Our algorithm {\tt AllGeConvexSets} renders all 200 subtrees of $T$ as a disjoint union of twenty two 012-rows. Generally the parameter $N=|GeConv(G)|$ in Theorem 1 can be replaced by the number $R$ of occuring 012-rows. While only the trivial bound $R\le N$ can be proven, in practise $R$ is often significantly smaller. The smaller $R$, the better the compression of $GeConv(G)$.

\vspace{3mm}
{\bf 4.2} If  only the $X\in r_2=(0,1,2,1,1,0,2,2,2,2)$ with $|X|=6$ are required, then we replace $r_2$ by $r_2':=(0,1,g_3,1,1,0,g_3,g_3,g_3,g_3)$. Generally the {\it g-wildcard}  $g_tg_t\cdots g_t$  signifies "exactly $t$ many 1's here". 

While the $g$-wildcard worked well in practice\footnote{If several $g$-wildcards occur within a single row, they must  be distinguished with extra indices. If however, as in our scenario, exactly one $g$-wildcard occurs per row, this $g$-wildcard may as well be dropped. For instance in $r_2$ it is clear that exactly $\binom{5}{3}$ bitstrings in $r_2$ have cardinality 6 (and they are easily pinpointed); furthermore exactly $\binom{5}{0}+\binom{5}{1}+\binom{5}{2}+\binom{5}{3}$ bitstings have cardinality $\le 6$.} in various scenarios, the {\it empty-row issue} usually prevents to prove the output-polynomial enumeration of all $k$-models. Namely, whenever a final 012-rows $r$ is "empty" in the sense of containing no $k$-models at all (because $r$ has too many 1's), $r$ is fruitless but generating it cost time. We return to that issue in 5.4.2.

\section{Compressing $Met(G)$ for arbitrary graphs $G$}

 After introductory remarks about Boolean functions (5.1), we show in 5.2 how the set of all geodesics of a graph gives rise to a (large yet useful) Boolean formula. We go on to illustrate the workings of {\tt AllMetricSets} on a lush toy example (5.3). Five loose ends, e.g. about the empty-row-issue, sampling at random, or chordless paths, are pressed into 5.4 and 5.5.

\vspace{2mm}
{\bf 5.1} Recall\footnote{For Section 5 the reader needs to know the rudiments of Boolean algebra as e.g. layed out in [M].} that a {\it clause} is a simple type of Boolean formula, i.e. a disjunction of literals like

$$(5)\quad \ol{x_1}\vee \ol{x_{12}}\vee x_2\vee x_3\vee x_7\vee x_6\vee x_5\vee x_{11}\vee x_{10}$$

If one replaces each {\it positive} literal $x_i$ in a clause by a {\it positive term}, i.e. a conjunction $x_ix_j\cdots x_k:=(x_i\wedge x_j\wedge\ldots\wedge x_k)$ of {\it positive} literals, one obtains a so-called {\it superclause} [W1,p.1083] like

$$(6)\quad \ol{x_1}\vee \ol{x_{12}}\vee x_2x_3x_7\vee x_2x_3x_8\vee\cdots\vee x_5x_{10}x_{11}$$

The superclausal algorithm  of [W2] is taylored to output,  by means of 012-rows, all models of a conjunction $\Psi$ of superclauses\footnote{While each superclause, as any Boolean formula, is equivalent [M,p.9]  to a conjunction of ordinary clauses, the latter may be exceedingly numerous.
.}. As we shall see, restriction to all models of Hamming-weight $\le k$ is possible as well.

\vspace{2mm}
{\bf 5.2}
Why is this algorithm relevant for us? Using depth-first-search it is shown in [W2]  that $G_1=(V_1,E_1)$ in Figure 1 has these  $s-t$ geodesics for $s=1,\ t=12$:

$$(7)\quad (s,2,3,7,t),\ (s,2,3,8,t),\ (s,2,6,7,t),\ (s,2,6,11,t),$$

$$(s,5,6,7,t),\ (s,5,6,11,t),\ (s,5,10,11,t)$$

\noindent
By definition of "metric" each $X\s V_1$ which contains $s,t$ and "wants" to be  metric, must 
contain one of the seven sets $\{2,3,7\},\{2,3,8\},...,\{5,10,11\}$. Observe that these sets match the positive terms in (6). Hence $X\s V_1$ passes the test with regards to the particular choice $(s,t)=(1,12)$ iff\footnote{More precisely, because (6) is also true when $\ol{x_1}\vee \ol{x_{12}}$ is true, $X$ passes the test also when $\{s,t\}\not\s X$.} its corresponding bitstring $x$ satisfies the superclause (6). 

If we set up  for any two non-adjacent vertices $s\neq t$ of $G_1$  a superclause $SC(s,t)$ analogous to (6), then $X\s V_1$ is metric iff $x$ satisfies the conjunction of all these formulas $SC(s,t)$. 

The bottom line is: Provided all geodesics are known [W2],  {\tt AllMetricSets}  reduces  to an instance of the superclausal algorithm. In Section 6 the role of Mathematica's command {\tt BooleanConvert} in all of this  comes to the fore.

\begin{figure}[htbp]
    \centering
    \includegraphics[scale=1.1]{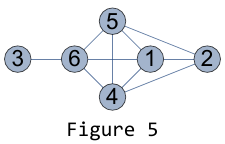}
\end{figure}

\vspace{2mm}
{\bf 5.3} To illustrate further, consider the graph $G_5$ in Figure 5.
We strive to enumerate all metric subsets $X\s V=\{1,2,...,6\}$.  Apart from the edges the five geodesics of $G_5$ induce this set $\Sigma$ of superclauses numbered 1' to 5':

\begin{itemize}
    \item[1' ] $\ol{x}_1\vee\ol{x}_3\vee x_6$
    \item[2' ] $\ol{x}_2\vee\ol{x}_3\vee x_1x_6\vee x_4x_6\vee x_5x_6$
\item[3' ] $\ol{x}_2\vee\ol{x}_6\vee x_1\vee x_4\vee x_5$
\item[4' ] $\ol{x}_3\vee\ol{x}_4\vee  x_6$
\item[5' ] $\ol{x}_3\vee\ol{x}_5\vee x_6$
\end{itemize}

\begin{tabular}{c|c|c|c|c|c|c|c}
	
       &1 & 2 & 3 &4 &   5 & 6 & comments  \\ \hline
	   & & & & & & &  \\ \hline\hline	
    
     $r=$ &2 & 2 & 2 & 2 & 2   & 2 & =powerset $\Pow(V)$ \\ \hline\hline
     
      $r_3=$ &{\bf 1} &  & {\bf 1}&  &    & 1 & {\bf final} of cardinality 8  \\ \hline
       $r_2=$ &{\bf 1} &  & {\bf 0} &  &    &  & {\bf final} of card=16 \\ \hline
      $r_1=$  &{\bf 0} &  & {\bf 2} &  &    &  & pending superclause 2'  \\ \hline\hline
      
      $r_1=$  &0 &  &  &  &    &  & p. superclause 2'  \\ \hline\hline

      $r_{13}=$  &0 & {\bf 1} & {\bf 1} &  &    &  & p. blackbox for $x_1x_6\vee x_4x_6\vee x_5x_6$  \\ \hline
      $r_{12}=$  &0 & {\bf 1}  &{\bf 0}  &  &    &  & p. superclause 3'  \\ \hline
       $r_{11}=$  &0 & {\bf 0} &{\bf 2}  &  &    &  & p. superclause 4'  \\ \hline\hline

       $r_{132}=$  &0 & 1 & 1 & 0 &  1  & 1 & {\bf final} of card=1 \\ \hline
        $r_{131}=$  &0 & 1 & 1 & 1 &  2  & 1 & {\bf final} of card=2 \\ \hline
      $r_{12}=$  &0 & 1  &0  &  &    &  & p.  superclause 3'  \\ \hline
       $r_{11}=$  &0 & 0 &2  &  &    &  & p. superclause 4'  \\ \hline\hline

        $r_{122}=$  &0 & 1  &0  &  &    & {\bf 1} & p. blackbox for $x_1\vee x_4\vee x_5$ \\ \hline
       $r_{121}=$  &0 & 1  &0  &  &    &{\bf 0}  & {\bf final} of card=4  \\ \hline
       $r_{11}=$  &0 & 0 &  &  &    &  & p.   superclause 4'  \\ \hline\hline

       $r_{1222}=$  &0 & 1 & 0 & 0 &  1  & 1 & {\bf final} of card=1 \\ \hline
        $r_{1221}=$  &0 & 1 & 0 & 1 &  2  & 1 & {\bf final} of card=2  \\ \hline
       $r_{11}=$  &0 & 0 &  &  &    &  & p.  superclause 4'  \\ \hline\hline

     $r_{11}=$  &0 & 0 &  &  &    &  & p.  superclause 4'  \\ \hline\hline

      $r_{113}=$  &0 & 0 &{\bf 1} &{\bf 1}  &    & 1 &   {\bf final} of card=2  \\ \hline
       $r_{112}=$  &0 & 0&  {\bf 1} & {\bf 0}  &    &  & p. superclause 5'  \\ \hline
        $r_{111}=$  &0 & 0 & {\bf 0} &{\bf 2}  &    &  & {\bf final} of card=8  \\ \hline\hline

         $r_{112}=$  &0 & 0 &1 &0  &    &  & p. superclause 5'  \\ \hline\hline

   $r_{1122}=$  &0 & 0 &1 &0  & {\bf 1}   & 1 & {\bf final} of card=1 \\ \hline
    $r_{1121}=$  &0 & 0&  1 &0  & {\bf 0}   &  & {\bf final} of card=2 \\ \hline\hline   
\end{tabular}

\vspace{2mm}
{\it Table 1: Snapshots of the LIFO stack underlying {\tt AllMetricSets}}
\vspace{2mm}

As to notation, considering say the 012-row $r_{131}$ in Table 1, we put $ones(r_{131}):=\{2,3,4,6\},\ zeros(r_{131}):=\{1\},\ twos(r_{131}):=\{5\}$. For better visualization we often replace all or some 2's by blanks; thus $twos(r_1)=\{2,3,4,5,6\}$.

We compute the modelset $Mod(\Sigma)$ by starting with the powerset $\Pow(V)=(2,2,2,2,2,2)=:r$
(see Table 1). It will shrink upon imposing the superclauses one by one. One option for superclause 1' to evaluate true is that either $x_1=0$ or ($x_1=1$ and $x_3=0$); this yields disjoint (012-)rows, i.e. $r_1\cup r_2=r_1\uplus r_2$. The second option $x_1=x_3=x_6=1$ is captured by row $r_3$. 

One checks that $r_2$ and $r_3$ happen to be {\it final} in the sense that they are wholly contained in $Mod(\Sigma)$. Upon removing and saving these (and all future final) rows, the LIFO stack consists of the single row $r_3$ which has superclause 2' pending to be imposed. As before, handling the negative literals is easy and yields $r_{11}\uplus r_{12}\s r_1$. If $x_2=x_3=1$ then $x_1x_6\vee x_4x_6\vee x_5x_6$ needs to be imposed on the row $(0,1,1,2,2,2)$. This task is handled by a "blackbox" (to be discussed later) which outputs $r_{131}\uplus r_{132}\s r_{13}$. 

And so it goes on. In the end we find that

$$Mod(\Sigma)=r_2\uplus r_3\uplus r_{132}\uplus\cdots
\uplus r_{1122} \uplus r_{1121}$$

\noindent
has $8+16+1+\cdots +1+2=47$ bitstrings; they bijectively match the metric subgraphs of $G_4$.

\vspace{9mm}
{\bf 5.4} Several remarks are in order. They concern the mentioned blackbox, the  enumeration of  cardinality-restricted models, distributed computation, and sampling at random.

\vspace{2mm}
{\bf 5.4.1} Recall that we invoked a "blackbox" to impose $x_1x_6\vee x_4x_6\vee x_5x_6$  upon the 012-row (0,1,1,2,2,2). Generally, imposing a disjunction of {\it positive} terms upon a 012-row boils down to the following problem. A given union $\rho_1\cup\cdots\cup \rho_t$ of 12-rows (which hence matches a set-filter generated by the $t$ sets $ones(\rho_i)$) must be partitioned as
$\rho_1\cup\cdots\cup \rho_t=\sigma_1\uplus\cdots\uplus \sigma_s$ (disjoint union). Here $s$ should be small and the $\sigma_i$'s are now 012-rows.

\vspace{2mm}
{\bf 5.4.2} As is well known, the essence of a LIFO stack (Last-In-First-Out) is that its occuring top rows repeatedly get replaced by $\ge 0$ many {\it candidate sons}. For instance in Table 1 the candidate sons of $r_1$ are $r_{11},r_{12},r_{131},r_{132}$. All of them happen to be {\it feasible} in the sense of containing models, two of them are even {\it final} (i.e. filled entirely with models). 

Conjunctions of superclauses being at least as difficult as conjunctions of clauses (=: CNF's), deciding the feasibility of a candidate son is NP-hard. But having a procedure that
answers "Is it feasibly?" by either "no" or "don't know", while avoiding false-negatives, is still helpful. 

Here it comes. Call a candidate son $k-${\it feasible} if it contains a model of Hamming-weight $\le k$. It then makes sense to throw away each candidate son $r_i$ of $r$ that has  $|ones(r_i)|>k$. Indeed,  the answer to "Is $r_i$ $k$-feasible?" is  "no" because the number of 1's can only increase upon further processing. (If it happens that {\it all} $r_i$'s get discarded, then $r$ itself was infeasible, albeit in hidden ways.)

The above mitigates the "empty-row issue" raised in 4.2 but doesn't sidestep it completely. A provingly output-polynomial enumeration of all $k$-feasible models is e.g. possible in [SW]. Finally observe that everything said about the Hamming-weight generalizes to arbitrary (positive) weight functions.

\vspace{2mm}
{\bf 5.4.3} As is well known, each depth-first search algorithm can (theoretically) be sped up  by any (w.l.o.g. integer) factor $fac$. This is best argued within the equivalent framework of LIFO stacks (see e.g. [W1,Sec.6.5]). In brief, when
 the initial  "head" computer has accumulated at least $fac$ items in its LIFO-stack, distribute them   in arbitrary manner to $fac$ many "satellite" computers. They launch their own LIFO-stacks, and in the end report back their results to the head.

\vspace{2mm}
{\bf 5.4.4} If either distributed computation cannot (by real-life obstacles) be launched to compute all metric sets, or the objective in the first place is to only  estimate $|Met(G)|$ (or other properties of $Met(G)$), then 
it is desirable to sample final 012-rows $r$ {\it uniformly at random}. As for any LIFO algorithm, this is easily achieved by  permuting the rows of the LIFO stack at random each time (or less frequently) its top row has been processed;
see also [W1,footnote 11].

\vspace{2mm}
{\bf 5.5} A path $(v_0,v_1,...,v_k)$ is {\it chordless} if $v_i,v_j$ are non-adjacent for all $0\le i<j\le k$. Clearly geodesic $\Ra$ chordless. It is easy to see that when the input for {\tt AllMetricSets} is inflated from the set of all geodesics to the set of all chordless paths, the output correspondingly inflates from $Met(G)$ to $Conn(G)$.

\section{Numerical experiments}

 This Section is about numerical experiments carried out with Mathematica; thus either {\tt AllMetricSets} or  {\tt AllGeConvSets} (abbreviated as "our algorithms" in the   heading of Table 2) are compaired with a method
   (abbreviated as "BConvert") which up to a few lines of extra code coincides with {\tt BooleanConvert}. Keep in mind, the latter being a hard-wired Mathematica-command, it has a head start on our "high-level" algorithms. 
   
  In 6.1 we further explain the notation of Table 2. In 6.2 we discuss the
 compression capabilities of the two competitors, and in 6.3 their CPU-times.
 Finally 6.4 is about {\tt SatisfiabilityCount}.

\vspace{2mm}
{\bf 6.1} Each row of Table 2 is dedicated to a single random graph.  
For instance $G=Gr(40,100,680)$ is a random graph with 40 vertices and 100 edges. It therefore triggers $\binom{40}{2}-100=680$ superclauses. Furthermore, $ m10: 15077(87s)$ signifies that it took {\tt AllMetricSets} 87 seconds to represent $Met(G,10):=\{X\in Met(G):\ |X|\le 10\}$ as a disjoint union of
 15077 many 012-rows $r_i$; so necessarily all $|ones(r_i)|\le 10$. Similarly $c10$ means that {\tt AllGeConvSets} was used instead.) 
The precise number of 
$k$-bounded metric subsets can be gleaned from the last column (see also 6.4). Thus for $Gr(40,100,680)$ and m10 the altogether 94394 metric subsets tell us that on average  each row  housed about 6 of them.

Recall that for trees "metric" is tantamount to "convex". For instance\\ $Tr(60,25,1711)$ is a random tree with 60 vertices (whence 59 edges) and 25 {\it leaves}. As before, $1711$ is obtained as $\binom{60}{2}-59$. It is obvious that a tree with $n$ nodes can have any number $x$ of leaves, where $x\in\{2,3,..,n-1\}$. As will be seen, the higher $x$ the higher the number of subtrees tends to be.

\vspace{2mm}
{\bf 6.2} Given any Boolean formula $fff$ the Mathematica command\\
{\tt BooleanConvert[fff,"ESOP"] } converts $fff$ into an {\it Exclusive Sum Of Products}, i.e. a DNF with $R$ terms $T_i$ such that  $T_i\wedge T_j=0$ for all $i\neq j$. In other words, we get a representation of the modelset $Mod(fff)$ as a disjoint union of $R$ many 012-rows. 

Our $fff$ of interest is the conjunction of the superclauses induced by a given random graph. Therefore {\tt BooleanConvert[fff,"ESOP"] } yields a representation of $Met(G)$ (or $GeConv(G)$) as disjoint union of 012-rows $r_i$. Scanning them all and  retaining exactly  the $r_i$'s with $|ones(r_i)|\le k$
yields such a representation for $Met(G,k)$. Unfortunately, in doing so one falls prey to the empty-row-issue (4.2). Surprisingly, this method (i.e. $BConvert$) nevertheless compressed somewhat better than our two algorithms. As to "somewhat", we mostly were within $200\%$ of $BConvert$ and the worst defeat\footnote{The original superclausal algorithm  (additionally to the don't-care "2") featured an extra $n$-wildcard which, for the particular (unbounded) instances in [W1], often beat {\tt BooleanConvert}.} occured for $Gr(60,1700,70)$ with 449'855 rows versus 74848 rows.

\vspace{2mm}
{\bf 6.3} What concerns CPU-times, the winner depends very much on $k$ and on the type of random graph evaluated. For instance the tree with 70 nodes and 30 leaves had its 9-bounded subtrees compressed by {\tt AllGeConvSets} in 58 sec, whereas $BConvert$ took 270 sec. Specifically, it took 81 sec to turn $fff$ into an ESOP with about 4.2 million terms (= disjoint 012-rows $r_i$). Afterwards it took 189 sec to pick the $r_i$'s with $|ones(r)|\le k$.
Things get worse (for both competitors) when increasing the number of leaves from 30 to 39.

However, for the majority of graphs in investigated\footnote{No attempt was made to render Table 2 statistically balanced in any way.} in Table 2  $BConvert$ wins out; including two cases where our algorithms could not terminate within 10 hours. But even in such situations many final 012-rows were delivered. In contrast, when $BConvert$ aborts, there are no partial results.

 \vspace{2mm}
 {\bf 6.3.1} Albeit (or because) $BConvert$ nicely compressed $Met(G)\ or\ GeConv(G)$ whenever it didn't abort,
 the author wondered whether one could  trim it further by mitigating the empty-row-issue. I pondered two approaches.

 The first  takes {\tt AllMetricSets} and replaces its blackbox subroutine (5.4.1)  by {\tt BooleanConvert[posterms,"ESOP"]}. Here
 {\tt posterms} is the conjunction of the positive terms of the  pending superclause. Unfortunately, by whatever reason, this hybrid algorithm was slower than the unhampered version of {\tt AllMetricSets}.

 The second approach relies on $bcf:=${\tt BooleanCountingFunction[n,k]} which is a very long\footnote{In fact $bcf$ has a whopping $\binom{n}{k}$ terms. To witness, if n=4 and k=2 then \\
$bcf=\ol{x}_1\ol{x}_2x_3x_4\vee\ol{x}_1\ol{x}_3x_2x_4\vee\ol{x}_1\ol{x}_4x_2x_3\vee\ol{x}_2\ol{x}_3x_1x_4\vee\ol{x}_2\ol{x}_4x_1x_3\vee\ol{x}_3\ol{x}_4x_1x_2$.} Boolean (ESOP) formula which evaluates to {\tt True} exactly for bitstrings of Hamming-weight $\le k$. It follows that
{\tt BooleanConvert[fff $\wedge$ bcf,"ESOP"]}  achieves what we want. Trouble is, due to the length of $bcf$ it was seldom efficient.

 \vspace{2mm}
 {\bf 6.4} The command {\tt SatisfiabilityCount} can calculate the {\it number} of (unboundend!) models of Gargantuan Boolean functions; thus it e.g. could excel for c80 in $Tr(80,35,3081)$ (so $k=80=n$). Trouble is, unless one relies on the cumbersome {\tt BooleanCountingFunction[n,k]} from above, {\tt SatisfiabilityCount} {\it cannot} count all $k$-bounded models.
 
 However, for most graphs appearing in Table 2 we are fine anyway. This is because from the ESOP delivered by $BConvert$ (or our algorithms) one gets the number of $k$-bounded models as a perk. Sometimes this ESOP was {\it soly} achieved by {\tt AllGeConvSets} (e.g. c10 in $Tr(80,35,3081)$), or {\it soly} by $BConvert$ (e.g. for $Gr(80,300,2860)$).
 
\vspace{1cm}

\begin{tabular}{|c|c|c|c|c|}
			
  graph or tree  &   {\tt our algorithms} & {\tt BConvert} & {\tt SatCount}   \\ \hline
                         &                        &                  &         \\ \hline
   
   $Tr(70,30,2346)$ & c9:$5999(58s)$ & $4436(270s)$  &  $94787$   \\ \hline
   $Tr(70,39,2346)$ & c9:$8946(353s)$ & $7811(1249s)$  &  $230'406$   \\ \hline
    $Tr(80,35,3081)$ & c80:$\ 5\cdot 10^6/row$ & aborted  & $6\cdot 10^{14}(35s)$ \\ \hline
$Tr(80,35,3081)$ & c10:$18071 (203s) $ & aborted  & 281'278 
\\ \hline\hline

$Gr(40,100,680)$ & m10:15077(87s)  &6790(1.4s)  &94394\\ \hline
$Gr(40,100,680)$ & c10:453(9.7s)  &316(0.1s)  &951\\ \hline
$Gr(40,200,580)$ & m10:247'043(966s)  &121'793(53s)  &1'329'829\\ \hline
 $Gr(40,200,580)$ &c10:397(9s)  &280(2s)  &634\\ \hline
 $Gr(80,300,2860)$ & m12: 22/row  &160217(8115s)  &3'262'493\\ \hline
$Gr(80,300,2860)$ & c5:1643(38s)  &1162(45s)  &3088\\ \hline \hline

$Gr(40,700,80)$ & m10:11462(7s)  &7669(0.4s)  &$10^9$\\ \hline
$Gr(40,700,80)$ & c10:65855(7.1s)  &35'080(2.3s)  &$14\cdot 10^6$\\ \hline
$Gr(60,1700,70)$ &m5:8976(5s)  &4506(0.2s)  &5'983'779\\ \hline
$Gr(60,1700,70)$ &m30:9003(5s)  &4506(0.3s)  &$6\cdot 10^{17}$\\ \hline
$Gr(60,1700,70)$ &c5:449'855(121s)  &74848(123s)  &4'054'335\\ \hline
$Gr(70,2360,55)$ &c5:236'953(63s)  &93'479(853s)  &10'446'417\\ \hline
\end{tabular}

\vspace{2mm}
{\sl Table 2: Numerical enumeration of $Met(G)$ and $GeConv(G)$}

\vspace{5mm}
\section{References}

\begin{itemize}
 \item[ {\bf AF} ] D. Avis, K. Fukuda, Reverse search for enumeration, Discrete Appl. Math. 65 (1996) 21-46.
 \item[ {\bf ASA} ] M. Alokshiya, S. Salem, F. Abed, A linear delay algorithm for
enumerating all connected induced subgraphs, BMC Bioinformatics, 20:319 (2019).
 \item[ {\bf BM} ] H.J. Bandelt, Mulder, Distance-hereditary graphs, Journal of Combinatorial Theory 41 (1986) 182-208.

    \item[ {\bf EJ} ]  P. Edelman, R. Jamison, The theory of convex geometries, Geometriae Dedicata 19 (1985) 247-270.
     
     \item[ {\bf FH} ] S.Fortunato, D. Hric, Community detection in networks: A user guide, PHYSICS REPORTS: REVIEW SECTION OF PHYSICS LETTERS 659 (2016): 1-44.

\item[ {\bf H} ] E. Howorka, A characterization of Ptolemaic graphs, Journal of Graph Theory 5 (1981) 323-331.
    
     \item[ {\bf M} ]  E. Mendelson, Boolean Algebra ans switching circuits, Schaum's Outline Series,
     McGraw-Hill 1970.
    
      \item[ {\bf S} ]  A. Schrijver, Combinatorial Optimization, Springer-Verlag Berlin, Heidelberg 2003.
       \item[ {\bf SW} ] Y. Semegni, M. Wild, Compression with wildcards: All k-models of a Binary Decision Diagram, arXiV:170308511v6 (2023), submitted.
    \item[ {\bf WX} ] S. Wang, C. Xiao, Novel Algorithms for Efficient Mining of Connected Induced Subgraphs of a Given Cardinality, arXiv:2112.07197.v3.

 \item[ {\bf W1} ]  M. Wild, Compression with wildcards: From CNFs to orthogonal DNFs by imposing the clauses one-by-one, The Computer Journal 65 (2022) 1073-1087.
 \item[ {\bf W2} ]  M. Wild, Enumerating all geodesics, arXiv:2409.16955v1.
 
\item[ {\bf W3} ]  M. Wild, Compression with wildcards: All induced subgraphs that are (respectively connected, chordal, bipartite, or forests, submitted.

 \item[ {\bf U} ]  T. Uno, Constant time enumeration by amortization, WADS 2015, LNCS 9214, pp. 593–605, 2015.
       
\end{itemize}

\end{document}